\documentclass[aps]{revtex4}
\usepackage{epsfig}
\usepackage{amsmath,amssymb,color}
\usepackage{epstopdf}
\parskip=\medskipamount

\newcommand{\be}{\begin{equation}}
\newcommand{\ee}{\end{equation}}

\newcommand\disp{\displaystyle}

\begin{document}

\title{Optimal control for virus spreading for an SIR model}

\author{O.D. Klimenkova$^{1}$}
\address{$^{1}$Department of Applied Mathematics, National Research University Higher School of
Economics, 101000, Moscow, Russia.}
%\date {November, 2016}
%\mail{odklimenkova@edu.hse.ru}

\begin{abstract}
In this paper, an SIR epidemic model with variable size of
population is considered.
We study optimal control problem for an SIR model with "vaccination" and "treatment" as controls.
It is shown that an optimal control exists.
We have already used functional, that lots of researchers use, and found that this functional is not appropriate, it has a defect.
Now, we fixed this defect by changing this functional.
We analyze the dependence of solutions on parameter of problems and discuss our result.
\end{abstract}

\maketitle

\section{Introduction}
One of the main method to investigate the process of virus infection in computer network is using mathematical epidemic models.
There are lots of epidemic models for
human disease, they can consider, for instance, incubation period,
the appearance of natural immunity,
natural mortality and birthrate of
individuals.
It should be mentioned that while using different epidemic models for computer network such things as the appearance of natural immunity and some other things are impossible.
Since mathematical models of virus spreading are described by system of nonlinear differential equations, lots of results can be obtained only numerically \cite{ZH}, \cite{kar}, \cite{zhang}.
But the most important part of modelling is not choice a numerical method or environment for calculation, it is a correct applying all parameters of model for particular problem.

We will consider an SIR model. We apply this model for computer networks. In this model nodes (computers) are divided by three groups: Susceptible, Infected and Recovered. Susceptible nodes can be infected by virus. Infected nodes are already infected. Removed nodes are nodes which are cured, for example antivirus was installed on this computers. This model can be described by following equations:
\begin{equation*}
\frac{dS}{dt} = - \beta SI,	
\frac{dI}{dt} = \beta SI - uI,	
\frac{dR}{dt} = uI
\end{equation*}
Here $u$  - intensity of "treatment" Infected nodes, $\beta$ - intensity of transmission of the virus for Infected to Susceptible nodes.

In this work we consider an SIR model with "vaccination" and "treatment". Supposed that except "treatment" of Infected nodes which means that viruses is removed and anti-virus program is installed, we have an ability to install anti-virus program on Susceptible nodes \cite{yusuf}.

The paper is organized as follows.
In section 2, we present an SIR model to be investigated and formulate an optimal control problem for that model.
In section 3, we derive the optimality system using Pontryagin's maximum principle and find structure  of optimal control.
In section 4, we solve the resulting optimality system numerically and
discuss our results.

\section{Model and optimal control problem}
We consider an SIR model with "vaccination" and "treatment" which can be described with the following differential equations:
\begin{equation}
\frac{dS}{dt} = - \beta SI - u_1 S,
\end{equation}
\begin{equation}
\frac{dI}{dt} = \beta SI - u_2 I - \alpha I,
\end{equation}
\begin{equation}
\frac{dR}{dt} = u_1 S + u_2 I
\end{equation}
Here $u_1$ is the proportion of the susceptible that is vaccinated per unit time, $u_2$ is the proportion of the infected that is treated per unit time, and $\alpha$ is the disease-induced death rate.
Initial conditions is:
\begin{equation}
S(0)=S_0, I(0)=I_0, R(0)=R_0
\end{equation}
We consider $u_1(t)$ and $u_2(t)$ as controls. Assume that admissible controls are measurable, bounded functions:
\begin{equation}
U=\{ (u_1 (t), u_2(t)):0\leqslant u_1(t) \leqslant u_{1max}, 0\leqslant u_2(t) \leqslant u_{2max} \}, t\in[0;T]
\end{equation}

Consider the following optimal control problem. In \cite{klim} we studied the problem of minimizing functional 

$\int_{0}^{T}(C_1 I(t) + C_2 u_2 ^2 + C_3 u_1 ^2)dt$ subject to (1)-(5).
We showed that functional that lots of researchers use \cite{yusuf}, \cite{bakare}  has a defect: when parameter $\alpha$ is growing up, this functional is going down, it means that high level of disease-induced death rate has a beneficial effect on functional.
For this reason that functional is not corresponding for us. Assuming that, the goal of this work is to introduce a new functional, without that defect.

To that end, we define a new group of nodes - D(t) (Defective) nodes which are defected by virus. Obviously:
$$\frac{dD}{dt}=\alpha I, D(0) = 0$$
Then we have:
$$D(T) = \int_0^T{\alpha I(t) dt}$$
Then we minimize the number of Defective, considering cost of "vaccination" and "treatment".
$$Z = \int_0^T (C_1 u_1^2 +C_2 u_2^2)dt + C_3 D(T) \longrightarrow \min $$
or
\begin{equation}
Z = \int_0^T (C_1 u_1^2 +C_2 u_2^2 + C_3 \alpha I(t))dt \longrightarrow \min
\end{equation}

$\textit{Remark.}$ Objective functional does not depend on variable R and equation (2) and (3) does not include this variable. So, we can reduce dimension of problem and excluded from consideration variable R and equation (3).

\subsection{Existence of solution}
$\textit{Theorem.}$ An optimal solution for the problem (1)-(6) exists.

Proof:

We apply Fillipov theorem \cite{Knowles}, we should show that:

1) set of acceptable solutions are bounded.

2) set of controls are convex compact.

3) velocity vector is convex by control.

Conditions 2 and 3 are held obviously. Let us show, that condition 1 is also held.

$\textit{Lemma.}$ Set of solutions of system (1)-(3) is bounded.

Proof:
From (1)-(3) we have $S(t) = S_0 e^{-\int_0^t (\beta I(\tau) + u_1 )d\tau } > 0$,
$I(t) = I_0 e^{-\int_0^t (\beta S(\tau) - \alpha - u_2 )d\tau } > 0$, $R(t)\ge 0$.
Also, we know, that $I(t) + S(t) + R(t) < N(t)$, consequently $S(t) < N(t)$, $I(t) < N(t)$, $R(t) < N(t)$ and $S$, $I$ and $R$ are bounded.
$\Box$
\subsection{Pontryagin maximum principle}
We apply Pontryagin Maximum Principle \cite{Pont} to the problem (1)-(6). Define Hamiltonian:
\begin{equation}
H=-\lambda_0 (C_1 u_1^2 +C_2 u_2^2 + C_3 \alpha I(t))+ \psi_1(t)(- \beta S(t)I(t) - u_1(t) S(t)) +\psi_2(t)(\beta S(t)I(t) - u_2(t) I(t) - \alpha I(t))
\end{equation}

Let $(S^*(t), I^*(t))$ - optimal solution in problem, $u_2^*(t), u_2^*(t)$ - corresponding optimal controls. Then according to Pontryagin Maximum Principle there will be found a constant $\lambda \ge 0$ and $\psi (t)=(\psi_1(t), \psi_2(t))$ such that:
\begin{enumerate}
	\item
$\frac{d\psi_1}{dt} = -\frac{\partial H}{\partial S} = \psi_1 \beta I + \psi_1 u_1 - \psi_2 \beta I$

$\frac{d\psi_2}{dt} = -\frac{\partial H}{\partial I} = \lambda_0 C_3\alpha +\psi_1 \beta S - \psi_2 \beta S + \psi_2 u_2 + \psi_2 \alpha$
	\item  $\psi_1(T) = 0, \psi_2(T) = 0 $
	\item \begin{multline*} H(S^*(t), I^*(t), \psi_1(t), \psi_2(t), \lambda_0,  u_1^*(t), u_2^*(t)) = \\
\\ \max_{0\leqslant u_1(t) \leqslant u_{1max}, 0\leqslant u_2(t) \leqslant u_{2max}} H(S^*(t), I^*(t), \psi_1(t), \psi_2(t), \lambda_0,  u_1(t), u_2(t))
\end{multline*}
\end{enumerate}
Remark, that for our problem $\lambda_0 \ne 0$. Below we will put $\lambda_0 = 1$.

\section{Analysis of maximum condition}
Consider maximum condition:
\begin{multline*}
\max_{  0\leqslant u_1(t) \leqslant u_{1max}, 0\leqslant u_2(t) \leqslant u_{2max}} -(C_1 u_1^2(t) +C_2 u_2^2(t) + C_3 \alpha I^*(t)) + \psi_1(t)(- \beta S^*(t)I^*(t) - u_1(t)S^*(t)) + \\
+ \psi_2(t)(\beta S(t)^*I^*(t) - u_2(t) I^*(t) - \alpha I^*(t)) =\\
= -(C_1 u_1^*(t)^2 +C_2 u_2^*(t)^2 + C_3 \alpha I^*(t)) + \psi_1(t)(- \beta S^*(t)I^*(t) - u_1(t)S^*(t)) + \\
+ \psi_2(t)(\beta S(t)^*I^*(t) - u_2^*(t) I^*(t) - \alpha I^*(t))
\end{multline*}

Write out items which have control $u_1$:
$$-C_1 u_1 ^2(t) - \psi_1(t) u_1(t) S(t) \longrightarrow \max_{0\leqslant u_1(t) \leqslant u_{1max}} $$
or
$$C_1 u_1 ^2(t)+ \psi_1(t) u_1(t) S(t) \longrightarrow \min_{0\leqslant u_1(t) \leqslant u_{1max}} $$

Define it as function $F(u_1) = C_1 u_1 ^2 + \psi_1 u_1 S$.

Analogically for $u_2$:
$$-C_2 u_2 ^2(t)- \psi_2(t) u_2(t) I(t) \longrightarrow \max_{0\leqslant u_2(t) \leqslant u_{2max}} $$
or
$$C_2 u_2 ^2(t)+ \psi_2(t) u_2(t) I(t) \longrightarrow \min_{0\leqslant u_2(t) \leqslant u_{2max}} $$
Define it as function $G(u_2) = C_2 u_2 ^2 + \psi_2 u_2 I$.

Functions $F(u_1)$ and $G(u_2)$ are convex, consequently minimum can be reached in
stationary point, where derivative is equals to 0, if this point is admissible. If this point is not admissible, minimum is achieved at $u_1=0$ or $u_1=u_{1\max}$ for $F(u_1)$ and at $u_2=0$ or $u_2=u_{2\max}$ for $G(u_2)$.

Find stationary point for $F(u_1)$ and $G(u_2)$.
$$F^\prime(u_1)=2C_1 u_1+\psi_1 S$$
$$F^\prime(u_1)=0 \Leftrightarrow u_1=-\frac{\psi_1 S}{2 C_1}$$
$$G^\prime(u_2)=2C_2 u_2+\psi_2 I$$
$$G^\prime(u_2)=0 \Leftrightarrow u_2=-\frac{\psi_2 I}{2 C_2}$$
Define:
$$u_1^b=-\frac{\psi_1 S}{2 C_1},
u_2^b=-\frac{\psi_2 I}{2 C_2}.$$
Here, we have the following structure of optimal control:
%\[
%u_2^* (t) =
%\begin{cases}
%u_2^b,  åñëè $0\leqslant u_2^b \leqslant u_{2\max}$ \\
%0, & åñëè $u_2^b < 0$ \\
%$u_{2\max }$, & åñëè $u_2^b > u_{2\max}$
%\end{cases}
%\]

\be
u_1^* (t) =
\left\{\begin{array}{rcll} \disp u_1^b, & \quad 0\leqslant u_1^b \leqslant u_{1\max}
\medskip \\ \disp 0, & \quad u_1^b < 0 \medskip \\
\disp u_{1\max }, & \quad u_1^b > u_{1\max} \medskip
\end{array} \right.
\label{eq:04}
\ee

\be
u_2^* (t) =
\left\{\begin{array}{rcll} \disp u_2^b, & \quad 0\leqslant u_2^b \leqslant u_{2\max}
\medskip \\ \disp 0, & \quad u_2^b < 0 \medskip \\
\disp u_{2\max }, & \quad u_2^b > u_{2\max} \medskip
\end{array} \right.
\label{eq:04}
\ee
\section{Numerical solutions}
System of equation of Pontryagin Maximum Principle is the following:
\be
\begin{cases}
\frac{dS}{dt} = - \beta SI - u_1 S
\\
\frac{dI}{dt} = \beta SI - u_2 I - \alpha I,
\\
\frac{d\psi_1}{dt} = -\frac{\partial H}{\partial S} = \psi_1 \beta I + \psi_1 u_1 - \psi_2 \beta I
\\
\frac{d\psi_2}{dt} = -\frac{\partial H}{\partial I} = \lambda_0 C_3\alpha + \psi_1 \beta S - \psi_2 \beta S + \psi_2 u_2 + \psi_2 \alpha
\end{cases}
\ee
We have a structure of optimal control (8)-(9), and boundary conditions:
$$
S(0)=S_0,  I(0)=I_0, \psi_1(T) = 0, \psi_2(T) = 0
$$
Thus, we turn our optimal control problem to boundary-value problem of  Pontryagin Maximum Principle for system of 4 equation the 1st order.
We use shooting method and Runge-Kutta 4th order procedure.
Then we consider dependence of optimal value of functional on parameters, especially on $\alpha$.
\begin{figure}[h!]
\center{\includegraphics[scale=0.5]{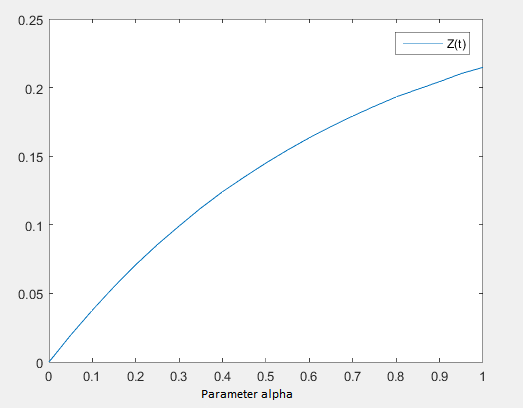}}
\caption{Dependence of functional Z on $\alpha$.}
\label{fig:NewF1}
\end{figure}

\section{Conclusion}
In this work the SIR model for virus spreading in computer networks is studied.
We consider optimal control problem for this model. We minimize a number of nodes which are defected by virus with "treatment" and "vaccination" as a controls. We prove existence of solution and find a structure of optimal control. We use shooting method and Runge-Kutta fourth order procedure for numerical solution. Then we analyze the dependence solution on parameter. We introduce functional without destructive dependence on disease-induced death rate.

\end{document}